\documentclass{article}
\usepackage{graphicx}
\usepackage[labelformat=empty]{caption}
\usepackage{amsmath}
\begin{document}
\begin{center}
\textbf{Cardinality Maximum Flow Network Interdiction Problem 
\\Vs. 
\\The Clique Problem}
\end{center}
\begin{center}
$Pawan Tamta^a,Bhagwati Prasad Pande^b,H.S.Dhami^c$
\end{center}
\begin{center}
a)Department of Mathematics, S.S.J Campus Almora,Kumaun University, Uttarakhand, India,pawantamta0@gmail.com.
\\b)Department of Information Technology, S.S.J Campus Almora,Kumaun University, Uttarakhand, India,bp.pande21@gmail.com.
\\c)Department of Mathematics, S.S.J Campus Almora,Kumaun University, Uttarakhand, India,profdhami@rediffmail.com.
\end{center}
\begin{center}
\textbf{Abstract}
\end{center}
Cardinality Maximum Flow Network Interdiction Problem (CMFNIP) is known to be strongly NP-hard problem in the literature. A particular case of CMFNIP has been shown to have reduction from clique problem. In the present work,an effort is being made to solve this particular case of CMFNIP in polynomial time. Direct implication of this solution is that the clique problem gets solved in polynomial time. 3-CNF Satisfiability and Vertex Cover problems, having reductions to and from the Clique Problem respectively, are also being solved in polynomial time by same algorithm. The obvious conclusion of the work is $P=NP$.
\begin{center}
\textbf{1. Introduction}
\end{center}
The maximum flow network interdiction problem (MFNIP) takes place on a network with a designated source node and a sink node. The objective is to choose a subset of arcs to delete, without exceeding the budget that minimizes the maximum flow that can be routed through the network induced on the remaining arcs.  
The study of MFNIP in particular originates from the Cold War. Now interdiction problems have many applications , including coordinating tactical air strikes [13], combating drug trafficking [16], controlling infections in a hospital [2], chemically treating raw sewage [14], and controlling floods [15]. 
From mid nineties to now, efforts have been made to develop some effective algorithms for MFNIP. Initially some naive algorithms were developed for interdiction problem such as a branch-and-bound strategy for general graph [7], and methods of varying quality for inhibition of s-t planar graph(planar graphs with both the source and sink on the outer face) [13]. 
Later in nineties efforts were made to categorize the problem and some polynomial time algorithms were developed on planar graphs for MFNIP. In 1993 Phillips [14] proved MFNIP as weakly NP Complete for planar graphs. At the same time Wood [16] introduced the Integer Linear Program (ILP) for MFNIP and proved it strongly NP Hard problem.
Recently Ricardo A. Collado et. al mentioned in Rutcor Research Report [6] that even the special case( of MFNIP),where the cost of arc removal is the same for each arc (CMFNIP) is known to be strongly NP-hard. It admits a very simple integer programming formulation [16]. A number of valid inequalities are known for this IP, but the integrality gap is still large [1]. The approximability of this problem is still unknown, with no positive or negative results in the literature. 
Rutcor Research Report [6] further envisages the recent results in the theory of Stackelberg games [3, 4, 12], which suggests that most of the network interdiction models are in fact APX-hard. Inapproximability bounds with a constant factor are known for shortest path interdiction problems,not for network flow interdiction problems.
In this paper we concentrate on Cardinality Maximum Flow Network Interdiction Problem (CMFNIP). CMFNIP is also known as k-most vital arc problem [15].CMFNIP is a special case of MFNIP with the restriction that interdiction cost for every arc is same [16]. Therefore in CMFNIP, we have to interdict the given number of arcs. CMFNIP is also known as strongly NP-hard problem [16]. 
We observe that if further restrictions are imposed on CMFNIP, then this problem can be solved in polynomial time. We name this problem as P-CMFNIP.
P-CMFNIP has been shown to have a reduction from Clique Problem [16, 1]. Therefore we get polynomial time solution for Clique Problem as well. Further 3-CNF Satisfiability and Vertex Cover problems have reductions to and from the Clique Problem [11] respectively. Therefore they are also being solved in polynomial time by same algorithm. We begin in section 2 with some preliminary definitions, the integer program given by Wood [16] and strengthened by Altner et. al. [1]. In section 3 the special case of CMFNIP named as P-CMFNIP is mentioned. An integer programming solution is proposed for P-CMFNIP. The integer program is then relaxed and shown to have zero integrality gap. In section 4 we mention the reduction of clique problem to P-CMFNIP. A polynomial time solution is given to the Clique Problem. Section 5 provides polynomial time solution to 3-CNF Satisfiability Problem and Vertex Cover Problem. Section 6 is about conclusion.    
\begin{center}
\textbf{2-Preliminaries}
\end{center}
A network is defined as $(N,A)$  where $N$ is the set of nodes and $A$ is the set of arcs. It is assumed that all of networks have a unique source $S\in N$ and a unique sink $t\in N$. Arc that originates from node $u$ and terminates at node $v$ are denoted by $(u,v)$. The $s-t$ cut is referred as  either a set of arcs that disconnects $S$ from $t$ upon their removal, or alternatively, as a bipartition of the nodes where $S$ and $t$ are not in the same partition. An undirected graph is denoted as $(V,E)$ where $V$ is the set of vertices and $E$ is the set of edges, an edge between vertices $u$ and $v$ by $\{u,v\}$ and an arc between node $i$ and $j$ as $(i,j)$. The capacity of every arc $(i,j)$ is denoted by $C_e$.The interdiction cost of any arc $e\in A$ is denoted by $r_e$ and total interdiction budget by $R$.
\\Wood [16] proposed the integer linear program for MFNIP and defined the decision variables as:
\begin{center}
\begin{equation*}
\alpha_v=
\begin{cases}
1&\text{if $v\in$$N$ is on sink side of the cut}\\
0&\text{otherwise}
\end{cases}
\end{equation*} 
\begin{equation*}
\beta_e=
\begin{cases}
1 &\text{if $e$$\in$$A$ is in the cut and is interdicted}\\
0&\text{otherwise}
\end{cases}
\end{equation*}
\begin{equation*}
\gamma_e=
\begin{cases}
1 &\text{if $e$$\in$$A$ is on sink side of the cut and is not interdicted}\\
0&\text{otherwise}
\end{cases}
\end{equation*} 
\end{center}
Integer linear program for complete formulation of MFNIP has been given by Wood[16] as under: 
\begin{equation}
\text{Minimize}\sum{C_e\gamma_e}\tag{2.1}
\end{equation}
\\subject to the conditions
\begin{equation}
\alpha_u-\alpha_v+\beta_(u,v)+\gamma_(u,v)\geq0\tag{2.2}
\end{equation}
\begin{equation}
\alpha_t-\alpha_s\geq1\tag{2.3}
\end{equation}
\begin{equation}
\sum_{e\in A} r_e\beta_e\leq R\tag{2.4}
\end{equation}
\begin{equation}
\alpha_v\in  \{0,1\} ,\forall v\in N\tag{2.5}
\end{equation}
\begin{equation}
\beta_e\in \{0,1\}, \forall e\in A\tag{2.6}
\end{equation}
\begin{equation}
\gamma_e\in \{0,1\},\forall e\in A\tag{2.7}
\end{equation}
Wood [16] further proved that MFNIP is a strongly NP-hard problem. He showed that CMFNIP (we call it P-CMFNIP) has a reduction from clique problem. Altner et. al. [1] obtained the natural linear programming relaxation for Wood's integer linear program[16] and showed that the integrality gap is not bounded below by any constant even when strengthened by valid inequality. 
\begin{figure}[ht]
\centering
\includegraphics[scale=0.5]{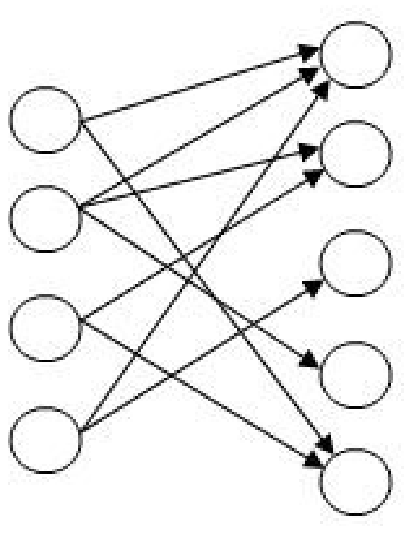}
\caption{figure 2.1}
\label{figure3}
\end{figure}
\\Altner et. al. [1] also proved that a simpler interdiction problem known as R-Interdiction Covering Problem (RIC) (figure 2.1) is strongly NP-hard.
Altner [1] et. al. showed that RIC has a reduction from Clique Problem. RIC has a very simple structure and it is identical to P-CMFNIP.
We observe that if CMFNIP will have finite number of node sets and one node set is connected to the next node set only ( as shown in figure 4.2), then it is solvable in polynomial time.
RIC and P-CMFNIP have two node sets only apart from source and sink node (figure 4.2).
Getting motivation from their special structure we modify the Integer Program of Wood [16] and Altner et. al. [1]. A zero integrality gap is obtained when we relax and strengthen this integer program.
\textbf{
\begin{center}
3- Linear programming solution to P-CMFNIP (Decision)
\end{center}
}
Referred to figure 4.2, we consider a particular case of CMFNIP and name it P-CMFNIP. P-CMFNIP was shown by Wood [16] to be a reduction from the Clique Problem. 
In this particular network we have a directed graph having a source node, a sink node and two node sets $A_1$ and $A_2$. Node set $A_1$ consists of the nodes connected to source node and similarly $A_2$ consists of the nodes directly connected to sink node. Furthermore, every node in $A_1$ is connected to exactly two nodes in $A_2$. The capacity of every arc connecting source node to $A_1$ has a capacity of 2 units, every arc connecting $A_1$ to $A_2$ has a capacity of 1 units, and every arc connecting $A_2$ to sink node has a capacity of 1. Interdiction cost of every arc is 1 so that for a set of arcs the interdiction cost is equal to the number of arcs in that set.
The decision version of P-CMFNIP is much simpler as compared to CMFNIP.  We are given constants $R$ and $K$. The decision version is given as; is it possible to interdict exactly $R$ nodes from node set $A_1$ to get the maximum flow in the remaining network equal to $K$.
\\Getting motivation from the structure of P-CMFNIP we modify the integer program of Wood[16].
\textbf{
\begin{center}
3.1-Formulation of Integer Program for P-CMFNIP
\end{center}
}
We define the decision variables as follws
\\Let $\gamma_i$  is the variable used for the $i_{th}$ node in $A_2$
\begin{equation*}
\gamma_i=
\begin{cases}
1&\text{if $i_{th}$ node does not gets removed}\\
0&\text{if $i_{th}$ node does gets removed}
\end{cases}
\end{equation*} 
Further let $\beta_j$ is the variable used for $j_{th}$ node in $A_1$
\begin{equation*}
\beta_j=
\begin{cases}
1&\text{if $j_{th}$ node is interdicted}\\
0&\text{if $j_{th}$ node is not interdicted}
\end{cases}
\end{equation*} 
We have to interdict nodes from $A_1$ only. The cut in the network is well defined. Therefore there is no need of $\alpha$ variables as defined by Wood[16].
\\Referred to the lemma 1 of section 3 by Wood[16] ; the maximum flow in P-CMFNIP is equal to the number of nodes in $A_2$, we design the objective function as under
\begin{equation}
\text{Minimize Z=}\sum{\gamma_i}\tag{3.1}
\end{equation}
Equation 3.1 is same as 2.1
\\We can interdict at most R nodes form node set $A_1$ where interdiction cost of every node is 1. The budget constraint is given as
\begin{equation}
\sum{\beta_j}\leq R\tag{3.2}
\end{equation}
Constraint 3.2 is same as 2.4 given by Wood[16].
\\We develop the third set of constraints by observing that any node in $A_2$ vanishes if all nodes in $A_1$ incident on it get interdicted, therefore we have
\begin{equation}
n\gamma_i + \sum_{j=1}^{n}\beta_j\geq n\tag{3.3}
\end{equation}
Constraint 3.3 is same as constraint 2.2. We have incorporated an additional condition in constraint 3.3. Here any node from $A_2$ can be interdicted by interdicting all nodes incident to it from $A_1$. Let $n$ nodes are incident on any $A_2$ node, then relation 3.3 is same as formulized in constraint 2.2 by Wood[16]. 
\\And the integer programming constraint is expressed as
\begin{equation}
\gamma_i, \beta_j \in \{0,1\} \forall i,j\tag{3.4}
\end{equation}
\textbf{
\begin{center}
3.2 Linear programming relaxation and strengthening
\end{center}
}
We relax the integer program by replacing constraint 3.4 by a weaker constraint
\begin{equation}
\gamma_i, \beta_j \in [0,1]\tag{3.5}
\end{equation}
According to the decision version of the problem, the minimum flow in the network remained after interdiction is at most $k$. Therefore simple constraint as given under is enough to strengthen the relaxed program.
\begin{equation}
\sum \gamma_i \geq k\tag{3.6}
\end{equation}
We denote the strengthened linear program given by equations \\$3.1$$,3.2$$,3.3$$,3.5$$,3.6$ as SLP.
\textbf{
\begin{center}
3.3 The strengthened linear program yields zero integrality gap
\end{center}
}
Let the decision problem (P-CMFNIP) has an affirmative answer and Z is the optimum solution obtained by the integer program (equations 3.1,3.2,3.3,3.4). Further let
\begin{equation}
Z=k\tag{3.7}
\end{equation}
We know that the nature of the objective function is to minimize therefore if $Z^*$ be the optimum solution obtained by SLP (equations 3.1,3.2,3.3,3.5,3.6), then it's obvious that
\begin{equation}
Z^*\leq Z=k\tag{3.8}
\end{equation}
Again by constraint 3.6 we have
\begin{equation}
Z^*\geq K\tag{3.9}
\end{equation}
From equations 3.7, 3.8 and 3.9 we conclude that 
\begin{equation}
Z=k=Z^*\tag{3.10}
\end{equation}
From equation 3.10 it's apparent that any optimum solution to the integer program given by equations 3.1, 3.2, 3.3, 3.4 will also be the solution to the SLP given by equations 3.1, 3.2, 3.3, 3.5, 3.6.
It is known that SLP being a simple linear program can be solved in polynomial time [9,10]. In next section we show that an integer solution to SLP can be decided in polynomial time. 
\textbf{
\begin{center}
3.4-Settlement of Integer solution to SLP
\end{center}
}
Let for some instance of PCMFNIP, the optimum solution obtained by the integer program (equations 3.1, 3.2, 3.3, 3.4) is $K$. Then the solution obtained by SLP (equations 3.1, 3.2, 3.3, 3.5, 3.6) is also $K$ (as proved in section 3.3). However one may argue that in case of alternate optimum solutions, a non integer solution can replace the integer solution. An integer solution can be decided in polynomial time. To support our assertion we prove a lemma.
\textbf{\\Lemma-} If the sum of the values of all variables involved yields an integer value then SLP can have an integer solution.
\\Proof- Let 
\begin{equation*}
\sum\gamma_r = K 
\end{equation*}
such that 
\begin{equation*}
\gamma_r \in (0,1) 
\end{equation*}
Where K is an integer.Further let $|\gamma_r|$=$K_1$ is greatest among all $\gamma$ variables. It is obvious that $K_1<1$ (as per supposition).
Next we find a value $K_2$ such that
\begin{equation*}
|\gamma_r| + K_2=1 
\end{equation*}
and we are left with $K-K_1-K_2$. We repeat the procedure by picking another largest non integer $\gamma$ variable unless we are left with $K-K_1-K_2-K_3\cdot\cdot\cdot\cdot=0$. All $\gamma$ variables that were picked in this manner are assigned value 1. Rest of the $\gamma$ variables remaining unpicked are assigned value 0. By this procedure sum of all $\gamma$ variables is still $K$ and they can have values either 0 or 1. Let the optimum solution consists $n$ $\gamma$ variables. Then it's clear that the whole procedure of rounding the variables to 0 or 1 is linear in $n$.
The same treatment can be given to $\beta$ variables also. In that case we replace $K$ by interdiction budget $R$. The worth mentioning fact is that the value of $R$ for P-CMFNIP is an integer. 
\textbf{
\begin{center}
3.5- Polynomial time solution to P-CMFNIP(Decision)
\end{center}
}
In this section we propose an algorithm solvable in polynomial time to solve PCMFNIP.
The algorithm is expressed as under
\textbf{\\Step1-} For a given directed graph H=(V,E), interdiction budget R and a constant K, SLP(equations 3.1, 3.2, 3.3, 3.5, 3.6) is solved.
\textbf{\\Step2-} If the optimum solution to SLP is K and the sum of variables involved is an integer then the given network has maximum flow K after interdiction. 
\\Step1 involves computation of the linear program (SLP) which can be solved in polynomial time[9,10]. Step2 being a simple if else statement can be decided in polynomial time. 
\begin{center}
\textbf{4-Reduction of the Clique Problem to P-CMFNIP}
\end{center}
In this section we reduce the Clique Problem  to P-CMFNIP, the reduction being given here is same as given in section 3 by Wood[16] .
The clique problem (decision) [16] is given as ; given an undirected graph H=(V,E)  and a positive constant K , does there exists a subgraph of H (complete graph) which is a clique on K vertices?
Here V is the set of nodes and E is the set of arcs. Clique of size K is a complete subgraph of H on K vertices i.e. $K\subset V$ such that every two nodes in it are connected by some arc in E.
For a given undirected graph H=(V,E) the reduction given by Wood[16] is as follows
\begin{figure}[ht]
\centering
\includegraphics[scale=0.6]{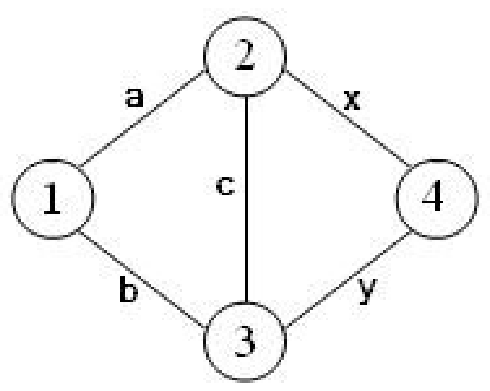}
\caption{figure 4.1}
\label{figure1}
\end{figure}
\begin{figure}[ht]
\centering
\includegraphics[scale=0.5]{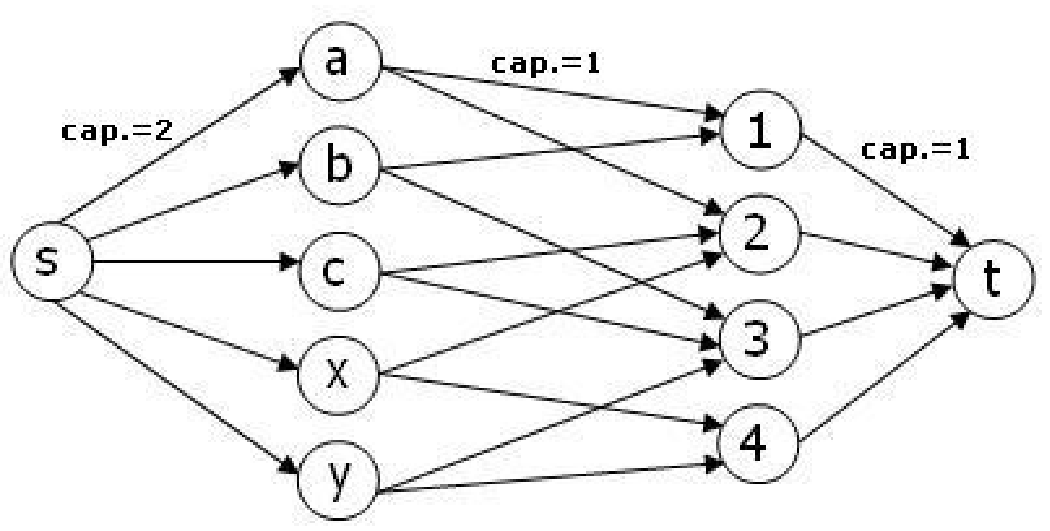}
\caption{figure 4.2}
\label{figure2}
\end{figure}
\\\\For each arc in E (figure 4.1) a node is constructed in $A_1$ (figure 4.2). Similarly for each node in V (figure 4.1) a node is constructed in $A_2$ (figure 4.2). Every node in $A_1$ is connected to exactly two nodes in $A_2$ (the idea is that one arc connects exactly two nodes). Every node in $A_1$ is connected to a source node and every node in $A_2$ is connected to a sink node. The interdiction cost of every arc is 1. The arc capacity of every arc connecting source node to $A_1$ has capacity 2, every arc connecting $A_1$ to $A_2$ has capacity 1, and every arc connecting $A_2$ to sink node has capacity 1.
Wood [16] showed that figure 4.1 contains a clique of size K if and only if the interdiction of $R= |E| - C_2^K$ nodes from $A_1$ yields the maximum possible flow of k units in the remaining network.
P-CMFNIP can be solved in polynomial time therefore any arbitrary instance of the clique problem can be solved in polynomial time.
\textbf{
\begin{center}
4.1- Polynomial time solution to the clique problem(Decision)
\end{center}
}
In this section we propose a polynomial time algorithm to solve clique problem (decision). 
The algorithm is based on linear programming formulation for clique problem. For that purpose we assign $\gamma$ variables to nodes. $\beta$ variables are assigned to arcs.
 Then the linear programming formulation for clique problem is given by SLP (equations 3.1, 3.2, 3.3, 3.5, 3.6). In constraint 3.3, $n$ is the degree of any $\gamma$ node. Summation is taken over all $\beta$ variables. These $\beta$ variables represent all those arcs which contribute to degree $n$.  
The algorithm is expressed as under
\textbf{\\Step1-} Given an undirected graph SLP is solved for given constant $k$ and $R= |E| - C_2^K$
\textbf{\\Step2-} If the optimum solution to SLP is K and the optimum solution can have the integer values of variables in 0 and 1 then the graph has a clique on K vertices, else it cannot have a clique on K vertices.
Both steps run in polynomial time as shown in section 3.5.
\textbf{
\begin{center}
4.2-Polynomial time solution to The Maximum Clique Problem (Optimization)
\end{center}
}
The optimization version of the Clique Problem is known as the Maximum Clique Problem. The problem is stated as; given an undirected graph H= (V,E), we have to find the complete subgraph of H of maximum size. Simply speaking we have to find a clique of maximum size.
The maximum number of vertices in a clique are $|E|$, therefore the algorithm runs as follows
\textbf{\\Step1-} SLP is computed for given value $K=|E|$ and interdiction budget $R= |E|- C_2^k$.  
\textbf{\\Step2-} If the solution is $K=|E|$ then the graph has a clique on $K=|E|$ vertices, else take $K=|E|-1$ and go to step 1.
\\A clique of size less than 2 in any undirected graph is not possible. Therefore in the loop of step 1 and step 2 the number of efforts cannot exceed $|E|-2$ which is polynomial in $|E|$. Step 1 and Step 2 run in polynomial time as shown in section 3.5.
\textbf{
\begin{center}
5-Polynomial time solution to 3-CNF Satisfiability Problem and Vertex Cover Problem
\end{center}
}
\begin{center}
\textbf{5.1 3-CNF Satisfiability Problem }
\end{center}
3-CNF satisfiability is defined by using the following terms. A literal in a boolean formula is an
occurrence of a variable or its negation. A boolean formula is in conjunctive normal form, or
CNF, if it is expressed as an AND of clauses, each of which is the OR of one or more literals.
A boolean formula is in 3-conjunctive normal form, or 3-CNF, if each clause has exactly
three distinct literals.
We define an instance of 3-CNF-SAT. 
Let $\phi = C_1\wedge C_2\wedge\cdot\cdot\cdot\cdot\wedge C_m$ be a boolean formula in 3-CNF with m clauses. For r = 1, 2,..., m, each clause $C_m$ has exactly
three distinct literals $l_1'$, $l_2'$, $l_3'$. AND operation is applied between any two clauses and OR operation is applied between two literals of any clause. We have to find the value of literals that makes the output of formula $\phi$ equal to 1.
\\It is obvious that K clauses contain 3K number of literals though some may be identical. Next we make the couple of every literal from each clause with every literal from the succeeding clause avoiding the coupling of any literal with its negation. The whole procedure runs in polynomial time as it is identical to the reduction used in[11]. The algorithm is based on linear programming formulation for 3-CNF Satisfiability problem. The linear programming formulation for 3-CNF Satisfiability Problem is given by SLP (Equations 3.1, 3.2, 3.3, 3.5, 3.6). For each couple we assign $\beta$ variable.$\gamma$ variables are assigned to literals. In third constraint (equation 3.3) $\gamma_j$ are all those variables which involve $\beta_i$ as member of couple and $n$ is the number of such $\beta$ variables.
The algorithm is given as under
\textbf{\\Step1- }SLP is computed for constant $K$ and $R= |E|-C_2^K$.
\textbf{\\Step 2- }  If the optimum solution to SLP is K and it can have integer solution, then we assign value 1 to these K variables and rest of $\gamma$ variables may have values 0 or 1. The formula thus obtained is satisfiable[11]. 
\\Both steps run in polynomial time as has been shown in section 3.5.
\textbf{
\begin{center}
5.2 Vertex Cover Problem
\end{center}
}
A vertex cover of an undirected graph $G = (V, E)$ is a subset $V'\subset V$ such that if $(u, v) \in E$,
then $u \in V'$ or $v \in V'$ (or both). That is, each vertex covers its incident edges, and a vertex
cover for $G$ is a set of vertices that covers all the edges in $E$. The size of a vertex cover is the
number of vertices in it.
The vertex-cover problem is to find a vertex cover of minimum size in a given graph.
Restating this optimization problem as a decision problem, we wish to determine whether a
graph has a vertex cover of given size k.
Vertex cover problem has reduction from clique problem[11]. Getting motivation from this concept we provide a direct polynomial time algorithm for vertex cover problem.
Given an undirected graph G=( V, E )  the complement of the graph is defined as G'=(V, E'). This means that any arc in E is not in E'.
As shown in [11] that G' has a clique of size V-K if and only if G has a vertex cover of size K.
\\The algorithm is based on linear programming formulation for vertex cover problem. For that purpose we assign $\beta$ variable to nodes and $\gamma$ variable to arcs belonging to G'.  
The linear programming formulation for Vertex Cover Problem is given by SLP (equations 3.1, 3.2, 3.3, 3.5, 3.6). In third constraint (equation 3.2)  $n$ stands for the degree of node $\gamma_i$ and variables $\beta_j$ stand for all those arcs in E' which contribute to degree $n$.

\noindent The algorithm is as given under
\textbf{\\Step 1-} SLP is computed for given constant $K$ and $R= |E|-C_2^{V-K}$.
\textbf{\\Step 2-} If the optimum solution is V-K and the optimum solution can have the integer values of variables in 0 or 1 then the given graph H=( V, E) has a vertex cover of size K.
\\As shown in section 3.5 both steps run in polynomial time. 
\textbf{
\begin{center}
6- Conclusion
\end{center}
}
The Clique Problem, Vertex Cover Problem and 3-CNF Satisfiability Problem are known to be NP-hard problems [11]. NP-hard problems can have polynomial time solution if and only if P=NP [5,11]. Based on the proof given in section 3 and algorithms proposed in section 3 and 4 we conclude that P=NP.
\textbf{\\References}
\\$[1]$ Douglas.S.Altner, Ozlem Eegun and Nelson A Uhan, The Maximum Flow Network Interdiction Problem, valid inequalities, integrality gaps and approximability, Operations Research Letter 38 (2010),pp 33-38.
\\$[2]$ N. Assimakopoulos, A network interdiction model for hospital infection control, Computers in Biology and Medicine 17 (1987), pp. 413-422. 
\\$[3]$ L. Bingol, A Lagrangian heuristic for solving a network interdiction problem, Master's thesis, Naval Postgraduate School, 2001.
\\$[4]$ Boyd, S. and Carr, R. (1999). A new bound for the ratio between the 2-matching problem and its linear programming relaxation. Mathematical Programming, Ser A, 86:499-514. 
\\$[5]$ Cook S.A. [1971], "The complexity of theorem-providing procedures", Proc. 3rdAnn. ACM symp. On Theory of Computing, Association for Computing Machinery, New York, 151-158
\\$[6]$ Ricardo.A. Collado , David Papp, Network Interdiction-Models, Applications, Unexplored Directions, Rutcor Research Report, RRR 4-2012, January 2012.
\\$[7]$ Ghare, P. M., Montgomery, D. C., and Turner, W. C. (1971). Optimal interdiction policy for a flow network. Naval Research Logistics Quarterly, 18:37-45.
\\$[8]$ Helmbold, R. L. (1971). A counter capacity network interdiction model. Technical Report R-611-PR, Rand Corporation, Santa Monica, CA.
\\$[9]$ Khachiyan, L.G., "Polynomial algorithm in linear programming,", Soviet Mathematics Doklady 20, (1979) pp. 191-194.
\\$[10]$ Karmarkar, N., "A new polynomial-time algorithm for linear programming," Combinatorica 4, (1984) pp. 373-395
\\$[11]$ Karp R.M. [1972], "Reducibility among combinatorial problems", in R.E Miller and J.W. Thatcher (eds.), Complexity of Computer Computations, Plenum Press, New York, 85-103.
\\$[12]$ Ford, L. R. and Fulkerson, D. R. (1962). Flows in Networks. Princeton University Press, Princeton, NJ.
\\$[13]$ A.W. McMasters and T.M. Mustin, Optimal interdiction of a supply network, Naval Research Logistics Quarterly 17 (1970), pp. 261-268.
\\$[14]$ C.A. Phillips, The network inhibition problem, in: Proceedings of the 25th Annual ACM Symposium on the Theory of Computing, 1993 pp. 776-785. 
\\$[15]$ H.D. Ratliff, G.T. Sicilia and S.H. Lubore, Finding the n most vital links in flow networks, Management Science 21 (1975), pp. 531-539.
\\$[16]$ R.K. Wood, Deterministic network interdiction, Mathematical and Computer Modelling 17 (1993), pp. 1-18. 
\end{document}